\documentclass[12pt]{amsart}
\usepackage{amsmath,amssymb,verbatim}
\newtheorem{thm}{Theorem}
\newtheorem{prop}{Proposition}[section]

\newtheorem{cor}[prop]{Corollary}
\newtheorem{defn}[prop]{Definition}
\theoremstyle{definition}
\newtheorem{exam}{Example}

\newcommand{\C}{{\mathcal C}}
\newcommand{\D}{{\mathcal D}}
\newcommand{\F}{{\mathcal F}}
\newcommand{\I}{{\mathcal I}}
\newcommand{\II}{\I^{(1)}}
\newcommand{\J}{{\mathcal J}}
\newcommand{\K}{{\mathcal K}}
\newcommand{\Lie}{{\mathcal L}}

\newcommand{\R}{{\mathbb R}}

\newcommand{\w}{\omega}
\newcommand{\st}{\widetilde{\sigma}}
\def\&{\wedge}
\newcommand{\di}{\partial}

\newcommand{\mood}{\ \operatorname{mod}\ }

\def\intprod{\mathbin{\raisebox{.4ex}{\hbox{\vrule height .5pt width 5pt depth 0pt %
        \vrule height 3pt width .5pt depth 0pt}}}}
\newcommand{\hook}{\intprod}
\newcommand{\bel}[1]{\begin{equation}\label{#1}}
\newcommand{\ee}{\end{equation}}

\begin{document}
\title{An Inverse Problem from Sub-Riemannian Geometry}
\author{Thomas A. Ivey}
\address{Dept. of Mathematics, College of Charleston,
Charleston SC 29424}
\email{ivey@math.cofc.edu}
\begin{abstract}The geodesics for a sub-Riemannian metric on a three-dimensional contact manifold
$M$ form a 1-parameter family of curves along each contact direction.  However, a collection of
such contact curves on $M$, locally equivalent to the solutions of a fourth-order ODE,
 are the geodesics of a sub-Riemannian metric only if a sequence of invariants vanish.
The first of these, which was earlier identified by Fels, determines if the differential equation is variational.
The next two determine if there is a well-defined metric on $M$ and if the given paths are its geodesics.
\end{abstract}
\dedicatory{To the memory of Robert B. Gardner}
\maketitle
\baselineskip=18pt
\section*{Introduction}
In this note, I will discuss the problem of recovering the geometric
structure of a three-dimensional contact manifold with a sub-Riemannian
metric from the geodesics for this metric.  (Sub-Riemannian metrics are also known
as Carnot-Carath\'eodory metrics.)  Since all the results herein
will be local in nature, the manifold may be taken to be an open
set $U\in \R^3$ with contact form $dy-z dx$, and we may assume that
on contact planes the metric has the form
$$g = E dx^2+F dx dz + G dz^2,$$
where $E,F,G$ are smooth functions on $U$ such that $g$ is positive definite.
The geodesics, as constructed via the Griffiths formalism,
 form a collection of paths tangent to the contact structure,
such that there is a 1-parameter family of distinct paths tangent to
each contact direction at each point.  Thus, part of the problem will be
to determine which such collections of paths come from a sub-Riemannian metric.

As explained below, the paths are locally equivalent to the integral curves of
a scalar fourth-order ODE.  The {\it variational multiplier problem} for
fourth-order ODE---i.e., the problem of characterizing equations
which are, up to multiple, the Euler-Lagrange equations for a second-order Lagrangian---was
solved by M. Fels \cite{fels}.  Since sub-Riemannian geodesics arise as solutions of
a variational problem, the present work is an extension of that of Fels;
to avoid confusion, the notation of \cite{fels} will be used whenever possible.

\section{Contact Path Geometries}
In this section I will review the construction of sub-Riemannian geodesics, and define
a $G$-structure canonically associated to the geodesics as paths.

Let $M$ be an oriented three-manifold with contact distribution $\D$ and
 sub-Riemannian metric $g$.  It is standard
 that one can associate to $g$
a $SO(2)$-structure $N$ inside the oriented coframe bundle $F(M)$, such
that, for any coframing which is a local section of $N$, the forms $(\w^1,\w^2,\w^3)$
of the coframing satisfy:
\begin{list}{(\roman{enumi})}{\usecounter{enumi}}
\item $\w^3$ annihilates the contact planes;
\item $(\w^1)^2+(\w^2)^2$ coincides with the metric on the contact planes; and,
\item $\w^1 \& \w^2$ gives the induced orientation.
\end{list}
Furthermore, one can choose $N$ uniquely so that there is
a connection form $\phi$ satisfying the structure
equations\footnote{This result appears in
\cite{ge}, where it is attributed independently to Bryant-Hsu and to G. Wilkens.}
\begin{align*}
d\w^1&= \phi \& \w^2 +(a_1 \w^1+a_2 \w^2) \& \w^3  \\
d\w^2&=-\phi \& \w^1 +(a_2 \w^1-a_1\w^2) \& \w^3  \\
d\w^3&= \w^1 \& \w^2 \\
d\phi&=K\w^1 \& \w^2 \mood \w^3.
\end{align*}
The functions $a_1,a_2$ are components of the torsion
of $g$ and $K$ is called its scalar curvature.\footnote{
It's clear that taking the $\w^i$ as an orthonormal coframe canonically associates
to $g$ with {\it Riemannian} metric $\widehat{g}$ on $M$, which induces $g$ on $\D$, and defines a
canonical foliation perpendicular to $\D$.  The torsion tensor
is the Lie derivative of $\widehat{g}$ along the leaves; if this vanishes, $g$ descends
to any (locally defined) quotient surface by foliation, and $K$ is
the Gauss curvature of the metric on that surface.}

Every contact curve in $M$ has a lift to $N$ on which the forms
$\w^2$ and $\w^3$ vanish.  (In this way, we'll identify $N$ with the space of
contact directions on $M$.)  Applying the Griffiths formalism \cite{griffiths}
to find the integral curves in $N$ of the Pfaffian system $\{\w^2,\w^3\}$ which are
extremal curves
for arclength $\int\w^1$, we obtain the following characterization of
sub-Riemannian geodesics\footnote{See \cite{monty} for earlier derivations of the geodesics by
other methods.}:
\begin{prop}\label{griffsub}  Let $Z$ be the rank one affine subbundle of $T^* N$
on which the canonical one-form is $\sigma=\w^1-x\w^3$,
$x\in\R$.  (Forms on $Z$ are pulled back via $\pi:Z\to N$.)
Then smooth geodesics are in 1-to-1 correspondence, via
the submersion $Z \to M$, with integral curves of the
Pfaffian system $\F=\{\theta_0,\theta_1,\theta_2,\theta_3\}$ on Z, where
\begin{align*}
\theta_0 &= \w^3\\
\theta_1 &= \w^2\\
\theta_2 &= \phi-x \w^1\\
\theta_3 &= dx-a_1 \w_1 -a_2 \w_2.
\end{align*}
\end{prop}

\noindent{\it Remark.}
In general, it is still an open question under what conditions all extremal curves for a given
variational problem with differential constraints arise as projections of integral curves of the differential system
formulated by Griffiths.  For example, in sub-Riemannian geometry in dimension four, exceptional
extremal curves exist which do not come from the Griffiths system.  Essentially, this is
because these {\it abnormal minimizers} \cite{monty} have few or no compactly supported
variations that are tangent to the given distribution.  However, by applying the regularity
test given by Hsu \cite{lucas}, one can show that for a sub-Riemannian metric on a contact
manifold, all geodesics arise via the Griffiths formalism.  (Intuitively, compact variations exist
because contact curves can be locally expressed in terms of an arbitrary function and its derivatives.)

\bigskip
Returning to the system $\F$ given above,
let $L$ be the line field on $Z$ which is annihilated by $\F$.
Integral curves of this line field
push down via $\pi$ to give a 1-parameter family of curves through each point of $N$, and push
down to $M$ to give a 1-parameter family of geodesics
tangent to each contact direction.
We will now generalize this situation, throwing away the metric:

\begin{defn}\label{bigdef} Let $M^3$ be a contact manifold.  Let $\rho:P \to M$ be a fibration,
with two-dimensional fibres, and $L$ be a line field on $P$ transverse to the two-dimensional fibres of $\rho$.
Let $\I$ be the Pfaffian system on $P$ which annihilates $L$ and the fibres of $\rho$,
and let $\J$ be the intersection of the retracting space \cite{EDS} of $\I$ with the annihilator of $L$.
Then $(P,L,\rho)$ defines a {\it contact path geometry} on $M$ if
\begin{list}{(\roman{enumi})}{\usecounter{enumi}}
\item for any vector $v \in L$, $\rho_*(v)$ is tangent to a contact plane;
\item the first derived system $\I'$ is one-dimensional at each point of $P$;
\item $\J'= \I$ at each point of $P$.
\end{list}
\end{defn}
\noindent
The last two conditions need explaining.  Because of the transversality of $L$, $\I$ is
two-dimensional.  Because of condition (i), $\I'$ contains the pullback of any contact form on $M$.  If $\I$
were integrable (i.e., $\I'=\I$, instead of being one-dimensional)
 then all paths through a given fibre $\rho^{-1}(x)$ would project down to a
single contact curve on $M$, so that there would be only one path through $x \in M$.
Condition (ii) implies that $\J$ is three-dimensional at each point.
It is automatic that $\I \subset \J'$.
If $\J$ were integrable, then integral surfaces of $\J$ would intersect $\rho^{-1}(x)$ in a 1-parameter family
of curves; since each such surface would project down to a single contact curve in $M$, this would imply
that there was only a 1-parameter family of paths through $x$.

Condition (ii) also implies that the three-dimensional distribution $\D$ containing $L$ and the kernel of $\rho_*$ is bracket-generating.
That in turn guarantees, by Chow's theorem \cite{chow}, that two arbitrary points in $M$ can be connected
by a piecewise smooth sequence of paths.

\begin{prop}
Given a contact path geometry we can construct, in a neighbourhood of any point $q \in P$, a coframe
 $(\sigma,\theta_0,\theta_1,\theta_2,\theta_3)$ such that
\begin{list}
{\arabic{enumi}.}{\usecounter{enumi}}
\item $v \in TP$ projects down to be a contact direction on $M$ if and only if $\theta_0(v)=0$
\item $\I=\{\theta_0,\theta_1\}$
\item \label{Jcond} $\J=\{\theta_0,\theta_1,\theta_2\}$
\item $L^\perp = \{\theta_0,\theta_1,\theta_2,\theta_3\}$
\end{list}
Moreover, these forms satisfy
\bel{gourst} d\theta_i \equiv \theta_{i+1} \& \sigma \mood \theta_0, \ldots, \theta_i, \qquad 0\le i\le 2.\end{equation}
These will be called {\it 0-adapted} coframes for
the contact path geometry.
\end{prop}
\begin{proof}
Let $\rho(q)=x\in M$.  On a neighbourhood $V$ of $x$, there exists a contact form  $\theta_0$, and 1-forms
$\sigma,\theta_1$ such that $d\theta_0\equiv \theta_1 \& \sigma \mood \theta_0$.
Pull these forms back to $U=\rho^{-1}(V)\subset P$; we will shrink $U$ when necessary.
Since $\sigma, \theta_1$ both restrict to be zero along the fibres of $\rho$, they cannot be independent
modulo $\I$.  Therefore we can arrange, by adding multiples of $\sigma$, that $\theta_1 \in \I$.
(Note that now $\theta_1$ is no longer the pullback of a form on $V$.).
Since $\theta_0,\theta_1 \in L^\perp$, then $\sigma \notin L^\perp$.
Since $\theta_0 \in \I'$, then $d\theta_1\ne 0 \mood \I$.

Since $\theta_0,\theta_1,\sigma$ span an integrable system, then there will be a smooth 1-form $\theta_2$
on $U$ that $d\theta_1 \equiv \theta_2 \& \sigma \mood \I$.  (Since $\I'$ is one-dimensional at each point,
$\theta_2$ is nonzero on $U$.)  By adding multiples of $\sigma$, we can arrange that $\theta_2 \in L^\perp$,
giving condition 3.
Because $\J' \ne \J$, there must be a nonzero 1-form $\theta_3$ on $U$ such that
$d\theta_2 \equiv \theta_3 \& \sigma \mood \J$.  We can similarly arrange that $\theta_3\in L^\perp$.
\end{proof}

\noindent{\it Remark.}
The above proposition could also be proved just using the assumption that $L$ is a line field
on $P$ and $\I,\J$ satisfy conditions (ii,iii) in Defn. \ref{bigdef}.  The contact structure and the submersion to $M$
can be recovered from $\I'$ and the retracting space $\C(\I')$ respectively.

\begin{cor}
In some neighbourhood $U$ of any given point $q \in P$, there exist coordinates
$x,y_0,y_1,y_2,y_3$ such that, for some function $F$ on $U$,
\begin{equation}\label{contactsys}
\begin{aligned}
\sigma &= -dx\\
\theta_0 &= dy_0-y_1 dx\\
\theta_1 &= dy_1-y_2 dx\\
\theta_2 &= dy_2-y_3 dx\\
\theta_3 &= dy_3-F(x,y_0,y_1,y_2,y_3) dx
\end{aligned}
\end{equation}
is a 0-adapted coframe.  Consequently, paths in $P$ are locally equivalent to the solutions of
the fourth-order ODE
\bel{fode}
y''''=F(x,y,y',y'',y''').
\end{equation}
\end{cor}
\begin{proof} The structure equations \eqref{gourst} enable us to apply the Goursat normal form
theorem \cite{EDS} to system $\J$.  This gives $\J=\{\theta_0,\theta_1,\theta_2\}$, in terms of the forms
defined here.  Since $\I=\{\theta_0,\theta_1\}$, then $dx \notin L^\perp$, and so there
exists some function $F$ such that $dy_3-F(x,y_0,y_1,y_2,y_3)dx=0$ along the paths in $U$.
\end{proof}

The set of 0-adapted coframes $(\sigma,\theta_0,\theta_1,\theta_2,\theta_3)$
for given contact path geometry forms a principal bundle over
$P$, with structure group $G_0 \subset GL(5,\R)$ consisting of matrices of the form
$$\begin{pmatrix} a&*&*&0&0\\0&b&0&0&0\\0&*&a^{-1}b&0&0\\
                  0&*&*&a^{-2}b&0\\0&*&*&*&a^{-3}b \end{pmatrix}.$$
This is precisely the $G$-structure that Fels associates to a
fourth-order ODE up to contact transformation (cf. \cite{fels}, Lemma 3.1).
Since the path geometry can be recovered uniquely from the $G$-structure, we will treat the
two notions as synonymous.

\section{Variational and Sub-Riemannian Path Geometries}
The goal of Cartan's method of equivalence \cite{RBG} is, for a given $G$-structure, to find a sub-bundle,
with reduced structure group, on which there exists a unique connection.  Like the Levi-Civita connection
in Riemannian geometry, this is typically obtained by fixing the value of all or part of the torsion
of the connection.
Then, invariants
may be extracted from the remaining torsion or the curvature of the connection.

We begin with Fels' result for $G_0$-structures of coframes satisfying \eqref{gourst}.  This gives a reduction of structure to the subgroup $G_1 \subset G_0$ consisting
of matrices of the form
$$\begin{pmatrix} a&0&0&0&0\\0&b&0&0&0\\0&0&a^{-1}b&0&0\\
                  0&0&0&a^{-2}b&0\\0&0&0&0&a^{-3}b \end{pmatrix}
\cdot \operatorname{exp}\begin{pmatrix} 0&0&0&0&0\\0&0&0&0&0\\0&c&0&0&0\\
                  0&0&\tfrac43c&0&0\\0&0&0&c&0 \end{pmatrix}.$$
In terms of path geometry, the result is:
\begin{thm}[Fels \cite{fels}]\label{felsthm}
Let $B_0 \searrow P$ define a contact path geometry.
Then there is a sub-bundle $B_1$
with three-dimensional structure group $G_1$, on which there exists a unique equivariant connection
satisfying the following structure equations:
\begin{subequations}\label{maineqs}
\begin{align}
d\sigma &= \alpha \&\sigma + \theta_0\&(T_1 \theta_1 + T_2 \theta_2 +T_3 \theta_3)
+ \theta_1 \& (T_4\theta_2 + T_5 \theta_3) \label{dsigma} \\
d\theta_0 &= \beta \&\theta_0 + \sigma \&\theta_1 \label{dthetazero} \\
d\theta_1 &= (\beta-\alpha)\&\theta_1 + \gamma \& \theta_0 +\sigma \&\theta_2 \label{dthetaone} \\
d\theta_2 &= (\beta-2\alpha)\& \theta_2+
\tfrac43\gamma\&\theta_1+ \sigma \&\theta_3 \label{dthetatwo} \\
d\theta_3 &= (\beta-3\alpha)\& \theta_3 +\gamma\&\theta_2 + \sigma \& (I_0\theta_0+I_1\theta_1)
+ T_6 \theta_0\&\theta_1 + T_7 \theta_0 \& \theta_2 + T_8 \theta_1 \& \theta_2. \label{dthetathree}
\end{align}
\end{subequations}
\noindent
[The one-forms $\alpha,\beta,\gamma$ are connection forms, and $I_0,I_1,T_1,\ldots T_8$ are
components of the torsion of the connection.]

Moreover, assuming $P$ is locally defined by a fourth-order ODE \eqref{fode},
 solutions of that ODE are
critical curves for a second-order Lagrangian if and only if the relative invariants $I_1$
and $T_5$ both vanish identically on $B_1$.  In that case, $T_8$ also vanishes.
\end{thm}
The essence of Fels' proof of the second statement is exhibiting a two-form on $B_1$,
$$\w = m(\theta_0 \& \theta_3 - \theta_1 \& \theta_2),$$
where $m$ is a non-zero function, such that $\w$ is closed and $G_1$-invariant.  (In fact,
$d\log m=3\alpha- 2\beta$, and, as Fels notes, the structure equations imply
that that one-form is closed in the variational case.)  It then follows
that $\w$ is the exterior derivative of the Poincar\'e-Cartan form associated to a Lagrangian
on the space of 2-jets.

We will speak of a path geometry for which $I_1,T_5,T_8$ vanish identically as
being {\it variational}.

\begin{exam}\label{ex1} Consider the second-order Lagrangian $\int e^{-3y''}dx$, for which
the Euler-Lagrangian equations are, up to multiple,
$$y''''-3(y''')^2=0.$$
The coframe \eqref{contactsys} gives a section of the bundle $B_0$
defining the corresponding $G_0$-structure on $J^3(\R,\R)$.
This coframe may be modified to give the following section of the reduced structure $B_1$:
\bel{lagframe}
\begin{aligned}
\theta_0 &= dy_0-y_1 dx\\
\theta_1 &= dy_1-y_2 dx\\
\theta_2 &= dy_2-y_3 dx-y_3\theta_1 +\tfrac3{10}y_3^2 \theta_0 \\
\theta_3 &= dy_3-3y_3 dy_2 -\tfrac3{10}y_3^2 \theta_1 +\tfrac65 y_3^3\theta_0\\
\sigma &= dx +\theta_1 -\tfrac35 y_3 \theta_0
\end{aligned}\end{equation}
Of course, the torsion satisfies $I_1=T_5=T_8=0$, but one also may compute\footnote{In
order to evaluate the torsion components along a given section of $B_1$, one must
determine the values of the connection forms in terms of the given coframe.  To do this, begin
with the $d\theta_0$ equation \eqref{dthetazero}, which determines $\beta$ modulo $\theta_0$.
One may set $\beta=\beta_0+b\theta_0$, where $\beta_0$ is any form satisfying \eqref{dthetazero}
and $b$ is not yet determined.  Then \eqref{dthetaone} determines $\alpha$ and $\gamma$
modulo $\theta_0,\theta_1$.  In fact, one may set
\begin{align*}\alpha&=\alpha_0+a\theta_0+z\theta_1\\
            \gamma &=\gamma_0-a\theta_1+c\theta_0.
\end{align*}
Now \eqref{dsigma} determines $z$ while \eqref{dthetatwo},\eqref{dthetathree} determine $a,b$ and $c$.}
that $T_2=\frac{12}5y_3$, $T_3=\frac35$ and $T_4=-1$ along this section of $B_1$.
\end{exam}

\begin{exam}\label{subrex}({\it sub-Riemannian geometry})
Let $Z$ be the five-manifold of Proposition \ref{griffsub}.
It is easy to verify that the 1-forms given there, when rounded out by $\sigma = \w^1-x\w^3$,
form a 0-adapted coframe for the corresponding contact path geometry.
We may adapt the coframe to obtain a section of the reduced bundle $B_1 \searrow Z$:
\bel{subrframe}\begin{aligned}
\theta_0 &= \w^3\\
\theta_1 &= \w^2\\
\theta_2 &= \phi-x\w^1+A \w^3\\
\theta_3 &= dx-a_1\w^1-(a_2+A)\w^2+B\w^3\\
\sigma &= \w^1-\frac35x\w^3
\end{aligned}\ \text{with}\ \begin{aligned}A &=\frac1{10}\left(a_2+3x^2-3K\right)\\
                                B &=\frac1{10}\left(s_2-3k_1-6a_1x-21b_1\right),
                    \end{aligned}
\end{equation}
where $K$ is the scalar curvature, and the $b_i$, $s_i$ and $k_i$ are defined on $N$ by
$$\left.\begin{aligned}
da_1 &\equiv 2a_2\phi +(s_1+b_2)\w^1+(s_2+b_1)\w^2\\
da_2 &\equiv -2a_1\phi +(s_2-b_1)\w^1+(b_2-s_1)\w^2\\
dK &\equiv k_1\w^1+k_2\w^2\end{aligned}\right\}\mood \w^3
$$
Again, one may compute that $I_1=T_5=T_8=0$, confirming that the path geometry is variational,
while $T_2=0$, $T_3=\frac35$, and $T_4=-1$ for this coframe.
\end{exam}
\bigskip
The fact that we obtained the same values for $T_3$ and $T_4$ as those from a general
 second-order Lagrangian hints at further relations among the torsion components.
One uncovers one of these while deriving the {\it refined structure equations}:

\begin{prop} Let $B_1$ be the canonical $G_1$-structure for a variational path geometry.
Then there exist functions $U_1,U_2$ on $B_1$ such that the connection forms satisfy
\bel{refalph}\begin{aligned}
d\alpha &= \tfrac23 d\beta\\
d\beta &= \sigma \& \gamma -\tau \& \theta_1-3\nu \& \theta_0 \\
d\gamma &\equiv \gamma \& \alpha -\tau \& \theta_2-\nu \& \theta_1 \mood \theta_0
\end{aligned}\ \text{where}\
\begin{aligned}
\tau&=T_1\theta_1+T_2\theta_2+T_3\theta_3\\
\nu &=U_1\theta_1+U_2\theta_2-T_2\theta_3+T_7\sigma
\end{aligned}.\end{equation}
The torsion components satisfy $T_3 =- \tfrac35 T_4$
and
\begin{align}
dT_1 &\equiv T_1(2\alpha-2\beta)-\tfrac43T_2\gamma - 2 U_1\sigma
& &\mood \theta_0,\theta_1,\theta_2,\theta_3 \\
dT_2 &\equiv T_2(3\alpha-2\beta)-\tfrac25T_4\gamma-(T_1+2U_2)\sigma
& &\mood \theta_0,\theta_1,\theta_2,\theta_3 \label{deeteetwo}\\
dT_4 &\equiv T_4(4\alpha-2\beta)-\tfrac53T_2\sigma
& &\mood \theta_0,\theta_1,\theta_2. \label{deeteethree}
\end{align}
\end{prop}
\noindent
The above equations indicate that $T_4$ is a {\it relative invariant} on $B_1$, i.e., it
varies along the fibres only by scaling.
Moreover, they indicate that the quadratic form $g=\sigma^2-T_4 \theta_1^2$ is well-defined,
up to multiple and modulo $\theta_0$, on $N$.  For, suppose $v$ is a vector field on $B_1$
which is annihilated by $\sigma,\theta_0,\theta_1,\theta_2$.  Then computing
the Lie derivative of $g$ gives
\begin{align*}
\Lie_v(g)&= 2\sigma \circ (v \hook d\sigma) -2T_4\theta_1 \circ (v \hook d\theta_1)-(v\hook dT_4)
\theta_1^2\\
    &\equiv 2(v \hook \alpha) \left[\sigma^2 -T_4 \theta_1^2\right]\mood \theta_0.
\end{align*}
This quadratic form will be our candidate for a sub-Riemannian metric.  Matters being so,
we will say that a variational geometry is {\it nondegenerate} if $T_4\ne 0$ everywhere,\footnote{Suppose
a variational path structure has $T_3$ and $T_4$ identically zero;
the refined structure equations show that $T_2=0$ also.
Recall that the system which restricts to be zero along the
fibres of $\rho:P\to M^3$ is spanned by $\sigma,\theta_0,\theta_1$.
Since $d\sigma \equiv T_1 \theta_0 \&\theta_1 \mood \sigma$,
vectors that are in the kernel of $\sigma$ push down to give a well-defined plane field on $M$.
These planes intersect the contact planes in a {\it distinguished} family of contact directions,
which are null lines with respect to $g$.}
and {\it definite} if $T_4$ is negative everywhere.
Assuming the latter is the case, then we may normalize $T_2$ and $T_4$ to have the
same values as in Example \ref{subrex}.

\begin{prop} Let $B_1\searrow P$ be a definite variational path geometry.  Then
there is a sub-bundle $B_2 \subset B_1$ on which
$$T_2=0 \ \text{and}\ T_4=-1.$$
On $B_2$ there exist smooth functions $W_0,W_1,W_2$, $G_0,G_1,G_2,G_3$, and $H$
such that
\begin{align}
\beta &=2\alpha + W_0\theta_0+W_1\theta_1+W_2\theta_2 \label{defw}\\
\gamma &= H\sigma - 3(G_0\theta_0+G_1\theta_1+G_2\theta_2+G_3\theta_3). \label{gamval}
\end{align}
\end{prop}
\begin{proof} Structure equations \eqref{deeteetwo},\eqref{deeteethree} show that
we may first pass to the sub-bundle where $T_4=-1$ and then move along
the fibres in a direction dual to $\gamma$ to pass to the sub-bundle where
$T_2=0$.  Once there, these equations show that $\beta-2\alpha$ and $\gamma$
restrict to have the above form.  Of course, \eqref{deeteetwo} shows that
$\frac25H=T_1+2U_2$.
\end{proof}

On $B_2$, the structure equations \eqref{maineqs} take the form
\bel{beeqns}\begin{aligned}
d\sigma &= \alpha \&\sigma + \theta_0\&(T_1 \theta_1 + \tfrac35 \theta_3)
- \theta_1 \& \theta_2\\
d\theta_0 &= 2\alpha \&\theta_0 + \sigma \&\theta_1\\
d\theta_1 &= (\beta-\alpha)\&\theta_1 + \gamma \& \theta_0 +\sigma \&\theta_2\\
d\theta_2 &= (\beta-2\alpha)\& \theta_2+
\tfrac43\gamma\&\theta_1+ \sigma \&\theta_3\\
d\theta_3 &= (\beta-3\alpha)\& \theta_3 +\gamma\&\theta_2 + I_0\sigma \& \theta_0
+ T_6 \theta_0\&\theta_1 + T_7 \theta_0 \& \theta_2.
\end{aligned}\end{equation}
with $\beta$ given by \eqref{defw}.

It's clear that the fibres of $B_2$ are one-dimensional, and the structure group
of $B_2$ is simply $\R^*$.
A element $\lambda \ne 0$ of this group acts on sections of $B_2$ by
$$g_\lambda \cdot (\sigma,\theta_0,\theta_1,\theta_2,\theta_3) =
(\lambda \sigma,\lambda^2\theta_0,\lambda\theta_1, \theta_2,\lambda^{-1}\theta_3).$$
Structure equations \eqref{refalph} show that $g_\lambda^*\alpha=\alpha$,
$g_\lambda^*\beta=\beta$, and $g_\lambda^*(\gamma) = \lambda^{-1}\gamma$.
Then the action on the new torsion is clearly
\begin{align*}g_\lambda \cdot (W_0,W_1,W_2) &= (\lambda^{-2}W_0, \lambda^{-1}W_1,W_2),\\
g_\lambda \cdot (H,G_0,G_1,G_2,G_3) &= (\lambda^{-2}H, \lambda^{-3}G_0,\lambda^{-2}G_1,\lambda^{-1}G_2,G_3).
\end{align*}
In particular, $W_2$, $G_3$, and the ratios $G_1:W_0$ and $G_2:W_1$ are invariant under the scaling action.

We should expect this scaling to be present, since two sub-Riemannian metrics
which differ by a constant factor have the same geodesics and hence define
the same path geometry.
For purposes of constructing a specific metric,
we will need to choose a section of $B_2$. Since
$3\alpha-2\beta$ is closed, integrals of this one-form comprise
a canonical codimension-one foliation of $B_2$ which is transverse to fibres and invariant
under the scaling action.

\begin{defn} A section of $B_2$ along which
\bel{candef}3\alpha-2\beta=0\end{equation}
will
be called a {\it canonical section} of $B_2$, or a {\it canonical coframe}
on $P$.  It follows from \eqref{defw} that
\bel{alfval}\alpha = -2(W_0 \theta_0+W_1\theta_1+W_2\theta_2)\end{equation}
along a canonical section.\end{defn}

One can check that the coframing constructed in Example \ref{subrex} is a canonical coframe.
Since such coframings are unique up to scale, it follows that if a path geometry
comes from a sub-Riemannian metric, then in terms of a canonical coframe
that metric must be $g=\sigma^2+(\theta_1)^2$.

\begin{prop}\label{metricprop}
Let $P$ be a definite variational path geometry on contact manifold $M^3$ and
$(\sigma,\theta_0,\theta_1,\theta_2,\theta_3)$ a canonical coframe on $P$.
Then $g=\sigma^2+(\theta_1)^2$ gives a well-defined metric on the
contact planes of $M$ if and only if $W_2$ is identically zero on $P$.
\end{prop}
\begin{proof} Let $v$ be any vector field on $P$ tangent to the fibres of the
projection $\rho:P \to M$.
Since $v$ is annihilated by $\theta_0,\theta_1$ and $\sigma$,
\begin{align*}
\Lie_v(g) &\equiv (v \hook \alpha) (\sigma)^2 + (v\hook (\beta-\alpha)) (\theta_1)^2 \mood \theta_0\\
    &\equiv -W_2 (v\hook \theta_2)\left( 2(\sigma)^2 + (\theta_1)^2 \right).
\end{align*}\end{proof}

Although the coframe \eqref{lagframe} is not a section of $B_2$, it can be adjusted so that $T_2=0$,
whereupon we see that $W_2$ is nonzero for Example \ref{ex1}.

For the rest of this section we will assume that $W_2$ is identically zero.
It remains to be seen if the given paths on $M^3$---which are projections
of the integral curves of the line field $L$---are the geodesics of the
sub-Riemannian metric we have constructed.  To investigate this further, we will need the
torsion identities
$$G_3=0, \quad G_2=W_1,$$
which result from computing $d(d\theta_1)=0$ using the structure
equations \eqref{beeqns}
and the equations \eqref{gamval},\eqref{candef},\eqref{alfval} with $W_2=0$.

\noindent{\it Remark.}\
One might wonder if other identities hold among the remaining torsion coefficients
$G_0$, $G_1$, $H$, $I_0$, $T_1$, $T_6$, $T_7$, $W_0$, $W_1$ as a result of our assumption that $W_2=0$.
However, no further identities arise, and this is proved by showing that the
exterior differential system defining a $G$-structure satisfying the structure
equations on $B_2$ with $W_2=0$ is {\em involutive}.

\begin{thm} Let $P$ carry a variational and definite path geometry with invariant
$W_2$ identically zero, and let $(\sigma,\theta_0,\theta_1,\theta_2,\theta_3)$
be a fixed canonical coframe on $P$.
Then the paths in $P$ project to be geodesics in $M$ for the sub-Riemannian
metric of Prop. \ref{metricprop} if and only if $G_1=2W_0$ identically on $P$.
\end{thm}
\begin{proof}
Let $N$ be the quotient of $P$ by the foliation by integral curves of the system
$\II = \{\sigma,\theta_0,\theta_1,\theta_2\}$.  Each contact curve in $M$
has a unique lift to $N$ as an integral curve of $\I=\{\theta_0,\theta_1\}$.  Clearly,
arclength is measured along these lifts by the integral of $\sigma$ modulo $\I$.  However, the
form $\sigma$ on $P$ does not descend to be well-defined on $N$, as shown by
\begin{align*}
d\sigma &= \alpha \&\sigma + \theta_0\&(T_1 \theta_1 + T_3 \theta_3)+T_4\theta_1 \& \theta_2\\
    &\equiv \tfrac35 \theta_0 \& \theta_3 \mood \Lambda^2 \II.
\end{align*}
However, computing $d^2\theta_0=0$ shows that
\bel{dwone}dW_1 \equiv (G_1 - W_0)\sigma + \tfrac15 \theta_3 \mood \I,\end{equation}
and this, together with $d\theta_0 \equiv 0\mood \Lambda^2\II$, shows that the 1-form
$$\st = \sigma + 3W_1 \theta_0$$
is well-defined on $N$.

Now arclength with respect to the metric may be measured on the integral curves of
$\I$ by the Lagrangian $\int \st$.  We will apply the Griffiths formalism \cite{griffiths} to investigate which
of these are geodesics for $g$.  Then, we will try to find conditions under which these curves coincide
with the projections of the paths in $P$ under $\pi:P\to N$.

Let $\xi = \st + x\theta_0 + y\theta_1$
on $Y=N \times \R^2$.  Then one finds that the two-form $d\xi$ is of full rank on $Y$, except where $y=0$.
Accordingly, let $\xi=\st + x\theta_0$ on $Z=N \times \R$.  Now one computes that
\bel{dixie}d\xi \equiv (dx + (3G_1-W_0) \sigma) \& \theta_0 +(\theta_2+(x+5W_1)\sigma)\&\theta_1
\mood \theta_0 \&\theta_1. \end{equation}
Let $\K$ be the
rank four Pfaffian system on $Z$ spanned by the four one-forms on the right in \eqref{dixie}:
$$\K=\{\theta_0,\theta_1,\theta_2+(x+5W_1)\sigma,dx + (3G_1-W_0) \sigma\}.$$

According to the Griffiths formalism, integral curves of $\K$ project to be extremal
curves for $\int \st$ on $N$.  These coincide with the projections of the paths
in $P$ if and only if, in a neighbourhood $U$ of each point of $P$, there is a local diffeomorphism
 $\varphi:U\to Z$ such that $\varphi^*\K$ coincides with $L^\perp=\{\theta_0,\theta_1,\theta_2,\theta_3\}$,
 the Pfaffian system on $P$ which defines the paths.  (The diffeomorphism would follow from the
 identification of paths with geodesics on $N$.)
The form $\varphi^*(\theta_2+(x+5W_1)\sigma)$ belongs in $L^\perp$ if and only if $\varphi^*x=-5W_1$.  Then, by
 \eqref{dwone},
$$\varphi^*(dx + (3G_1-W_0)\sigma) \equiv (4W_0-2G_1)\sigma \mood L^\perp,$$
showing that $\varphi^*\K = L^\perp$ if and only if $G_1=2W_0$.
\end{proof}

\section{Discussion}
The results of the previous section may be surprising.  For, one could reason that, once
a path geometry is known to be variational, it must arise from a second-order
Lagrangian, of the form
$$\int L(x,y,y',y'') dx,$$
satisfying the nondegeneracy condition $\di^2L/\di(y'')^2\ne 0$.  Then $L$ is of the
form
$$L(x,y,y',y'')=\sqrt{E + F y'' + G (y'')^2},$$
for some functions $E,F,G$ of $x,y,y'$, if and only if $L$ satisfies a certain third-order ODE
as a function of $y''$.  In other words, it seems like only one extra condition must be satisfied
in order for the Lagrangian to be length with respect to a sub-Riemannian metric.
Instead, we find
that two scalar conditions (in addition to the Fels variational condition) must hold in
order for the metric to be well-defined and in order for its extremals to coincide with
the given paths.  (The reader should note that the above remark about involutivity
implies that the condition $G_1=2W_0$ is independent from $W_2=0$.)
It would be interesting to find examples of variational
path geometries (equivalently, variational fourth-order ODE) for which $W_2=0$
but the extremals of the associated metric do not coincide with the given paths.  Such examples must exist,
again, because of the involutivity of $W_2=0$.
\section*{Acknowledgements}
I am grateful to Ian Anderson and Mark Fels for comments and suggestions.  Most of the calculations in this paper
have been made using Maple V and Maple 6.

\end{document}